\documentclass[a4paper,12pt]{article}
\usepackage[top=3cm,bottom=3cm,left=2.5cm,right=2.5cm]{geometry}
\usepackage{t1enc}
\usepackage[latin1]{inputenc}
\usepackage[english]{babel}
\usepackage{amsmath,amsthm}
\usepackage{amsfonts}
\usepackage{latexsym}
\usepackage[dvips]{graphicx}
\usepackage{graphicx}
\usepackage{color}
\usepackage{amsmath}
\usepackage{amssymb}
\usepackage{amsbsy}
\usepackage{amsfonts}
\usepackage{amsfonts,mathrsfs}

\begin{document}
\newtheorem{prop}{Proposition}[section]
\renewcommand{\theprop}{\arabic{section}.\arabic{prop}}
\newtheorem{lem}[prop]{Lemma}
\renewcommand{\thelem}{\arabic{section}.\arabic{lemma}}
\newtheorem{thm}[prop]{Theorem}
\renewcommand{\thethm}{\arabic{section}.\arabic{thm}}
\newtheorem{cor}[prop]{Corollary}
\renewcommand{\thecor}{\arabic{section}.\arabic{Cor}}
\newtheorem{conj}[prop]{Conjecture}
\renewcommand{\theconj}{\arabic{section}.\arabic{conj}}
\renewcommand{\theequation}{\arabic{section}.\arabic{equation}}
\renewcommand{\thefigure}{\arabic{figure}}
\newtheorem{exa}{prop}
\renewcommand{\theexa}{\arabic{section}.\arabic{ex}}
\newtheorem{rem}[prop]{Remark}
\renewcommand{\therem}{\arabic{section}.\arabic{rem}}
\renewcommand{\theequation}{
\arabic{equation}}

\newtheorem{defi}[prop]{Definition}
\renewcommand{\thedefi}{\arabic{section}.\arabic{defi}}
\renewcommand{\thefigure}{\arabic{figure}}
\def\zm{\noindent{\bf Proof.\ }}
\def\ezm{\vspace*{6mm}\framebox{} }

\addtolength{\baselineskip}{+0.6mm}


\title{Molecular trees with extremal values of the second Sombor index \thanks{This work is supported by the National Natural Science Foundation of China (No. 11971164) and  the Hunan Provincial Natural Science Foundation of China (2020JJ4423).}}
\author{Zikai Tang, Hanyuan Deng \footnote{Corresponding author: hydeng@hunnu.edu.cn}\\
{\small MOE-LCSM, School of Mathematics and Statistics, Hunan Normal University,}
 \\{\small  Changsha, Hunan 410081, P. R. China.}
}

\date{}
\maketitle

\begin{abstract}
A new geometric background of graph invariants was introduced by Gutman, of which the simplest is the second Sombor index $SO_2$, defined as $SO_2=SO_2(G)=\sum_{uv\in E}\frac{|d^2_G(u)-d^2_G(v)|}{d^2_G(u)+d^2_G(v)}$, where $G = (V, E)$ is a simple graph and $d_G(v)$ denotes
the degree of $v$ in $G$. In this paper, the chemical applicability of the second Sombor index is investigated and it is shown that the the second Sombor index is useful in predicting physicochemical properties with high accuracy compared to some well-established and often used indices. Also, we obtain a bound for the second Sombor index among all (molecular) trees with fixed numbers of vertices, and characterize those molecular trees achieving the extremal value.

{\bf Keywords}: Second Sombor index; Tree; Molecular tree; Extremal value.
\end{abstract}

\maketitle

\makeatletter
\renewcommand\@makefnmark%
{\mbox{\textsuperscript{\normalfont\@thefnmark)}}}
\makeatother

\baselineskip=0.25in

\section{Introduction}
Let $G$ be a simple connected graph with vertex set $V(G)$ and edge set $E(G)$, and denote by $n=|V(G)|$ and $m=|E(G)|$ the number of vertices and edges, respectively. The degree of a vertex $v$ in $G$, denoted by $d_G(v)$ or $d(u)$, is the number of its neighbors. If the vertices $u$ and $v$ are adjacent, then the edge connecting them is labeled by $e=uv$. In the mathematical and chemical literature, several dozens of vertex-degree-based (VDB) graph invariants (usually referred to as "topological indices") have been introduced and extensively studied \cite{gut-1,gut-2,gut-3,Ku20,TC09}. Their general formula is
\begin{equation}\begin{split} TI(G)=\sum_{uv\in E}F(d_G(u),d_G(v)),\nonumber\end{split}\end{equation}
where $F(x, y)$ is a function with the property $F(x, y) = F(y, x)$.

Recently in \cite{gut-2} Gutman showed that geometry-based reasonings reveal the geometric background of several classical topological indices (Zagreb, Albertson) and introduced a series of new Sombor-index-like VDB invariants, denoted below by $SO_1,SO_2, \dots, SO_6$. The second Sombor index $SO_2$ of a graph $G$ is defined as
\begin{equation}\label{eq-0}
\begin{split}
SO_2=SO_2(G)=\sum_{uv\in E}\frac{|d^2_G(u)-d^2_G(v)|}{d^2_G(u)+d^2_G(v)}.
\end{split}
\end{equation}

A molecular tree is a tree of maximum degree at most four. In this paper, the chemical applicability of the second Sombor index is investigated and it is shown that the second Sombor index is useful in predicting physicochemical properties with high accuracy compared to some well-established and often used indices. Also, a bound for the second Sombor index among all (molecular) trees with fixed numbers of vertices is obtained, and those molecular trees achieving the extremal value are characterized.

\section{The chemical applicability of the second Sombor index}

In this section, the chemical applicability of the second Sombor index is investigated. We consider the data set of octane isomers for such testing and corresponding experimental values of physico-chemical properties are collected from http://www.moleculardescriptors. eu/dataset/dataset.htm.
First, we give experimental values of the second Sombor index of octane isomers, which are listed in Table \ref{tb-1}, where there are two pairs of octane isomers with identical values of the second Sombor index since they have the same degree coordinates.

\begin{table}
\begin{center}
	\caption{ Experimental values of the second Sombor index of  for octane isomers} \label{tb-1}
	\begin{tabular}{cc|cc}
		\hline
Molecule& $SO_2$ &Molecule& $SO_2$  \\ \hline
octane&	1.2&
2-methyl-heptane&	2.5846\\
3-methyl-heptane&	2.7692&
4-methyl-heptane&	2.7692\\
3-ethyl-hexane&	2.9538&
2,2-dimethyl-hexane&	3.8471\\
2,3-dimethyl-hexane&	3.3846&
2,4-dimethyl-hexane&	4.1538\\
2,5-dimethyl-hexane&	3.9692&
3,3-dimethyl-hexane&	4.1647\\
3,4-dimethyl-hexane&	3.5692&
2-methyl-3-ethyl-pentane&	3.5692\\
3-methyl-3-ethyl-pentane&	4.4824&
2,2,3-trimethyl-pentane&	4.7117\\
2,2,4-trimethyl-pentane&	5.2317&
2,3,3-trimethyl-pentane&	4.8447\\
2,3,4-trimethyl-pentane&	4&
2,2,3,3-tetramethylbutane&	5.2941\\
		 		\hline
	\end{tabular}
\end{center}
\end{table}

By the experimental values of AcenFac, S, SNar and HNar of octane isomers (from http://www.moleculardescriptors. eu/dataset/dataset.htm.) and Table \ref{tb-1},  we find the correlation of AcenFac, S, SNar and HNar with the second Sombor index $SO_2$ for octane isomers. The data related to octanes are listed in Table \ref{tb-2}. The following equations give the regression models for the second Sombor index.
 \begin{equation}
 \begin{split}
 AcenFac=0.4536-0.0314\times SO_2
 \end{split}
 \end{equation}
 \begin{equation}
 \begin{split}S=119.1755-3.6697\times SO_2
 \end{split}
 \end{equation}
 \begin{equation}
 \begin{split}
 SNar=4.6576-0.3003\times SO_2
 \end{split}
 \end{equation}
 \begin{equation}
 \begin{split}
 HNar=1.7137-0.0815\times SO_2
 \end{split}
 \end{equation}

\begin{table}
\begin{center}
	\caption{ The square of correlation coefficient of the second Sombor index with AcenFac, S, SNar and HNar} \label{tb-2}
	\begin{tabular}{c|cccc}
		\hline
&AcenFac&S& SNar & HNar\\ \hline
$SO_2$&0.9202&0.8433&0.9356&0.9512\\
\hline
	\end{tabular}
\end{center}
\end{table}

\begin{figure}[th]
	\centering
\includegraphics[width=0.48\textwidth]{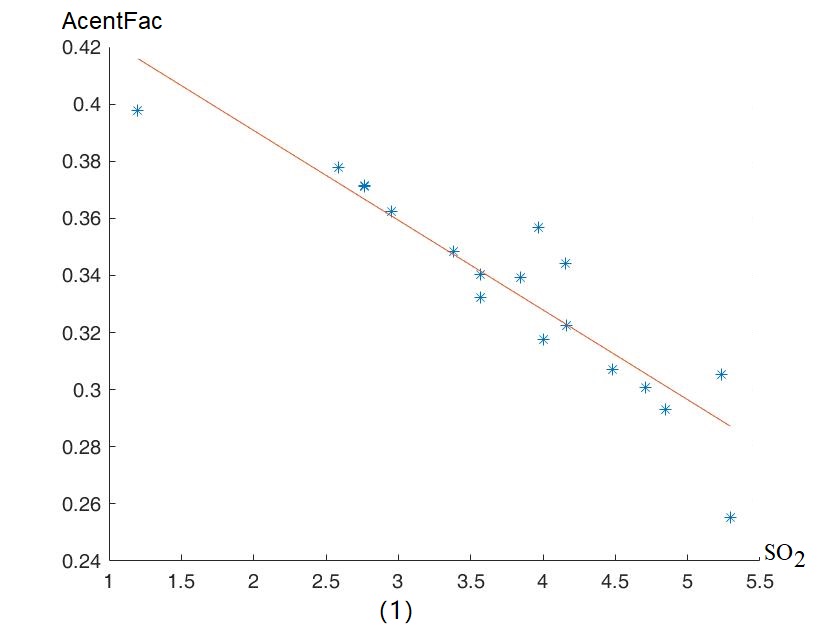}
\includegraphics[width=0.48\textwidth]{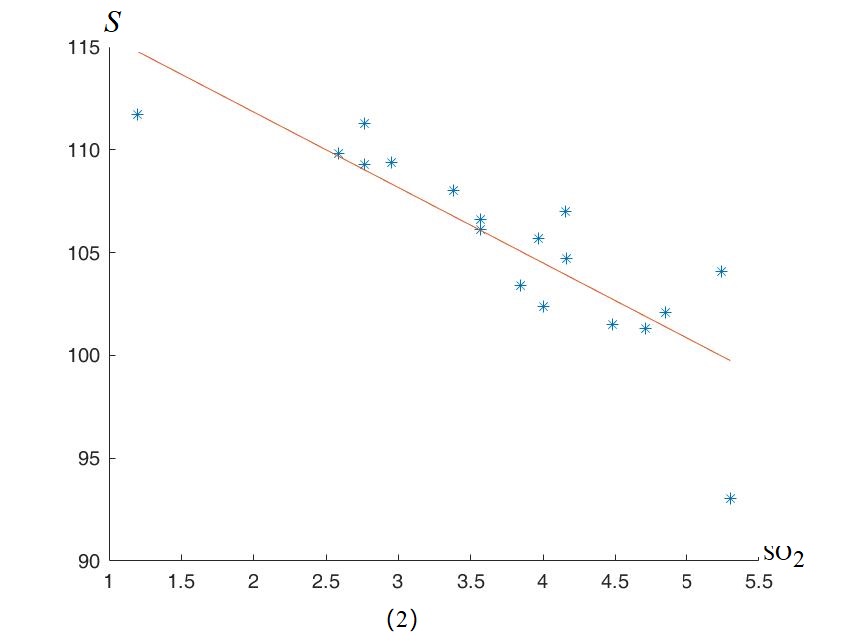}\\
\includegraphics[width=0.48\textwidth]{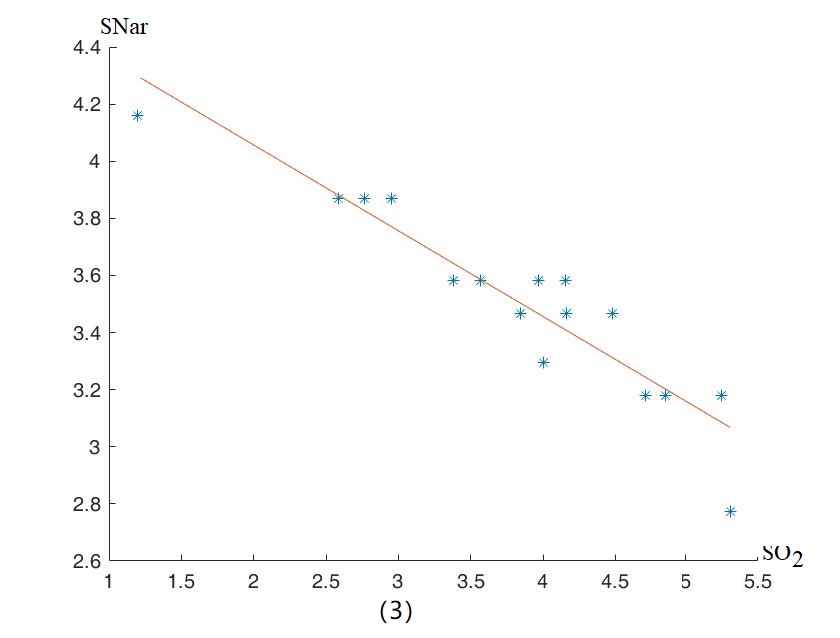}
\includegraphics[width=0.48\textwidth]{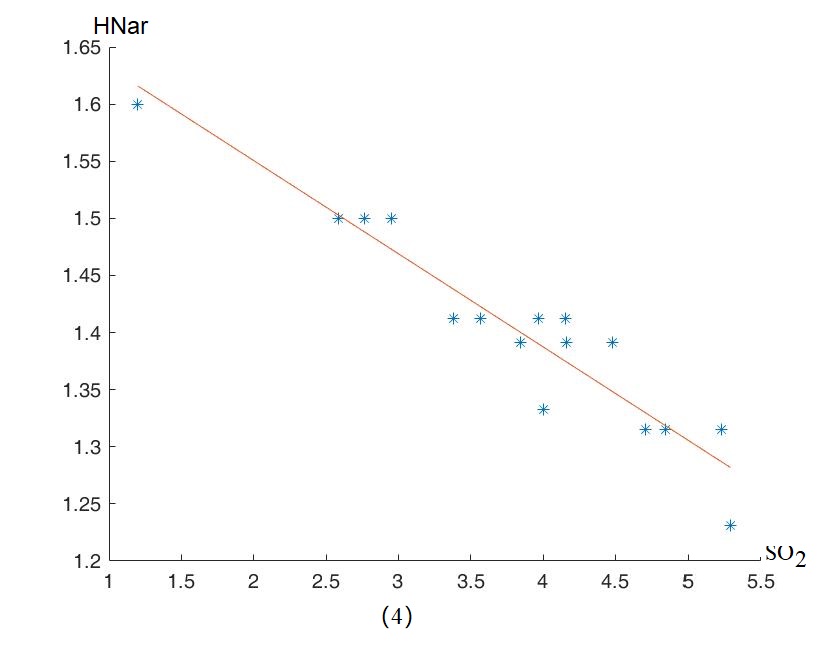}

	\caption{Linear fittings of the second Sombor index ($SO_2$) with AcenFac, S, SNar and HNar for octane isomers}\label{fig-1}
\end{figure}

Figure \ref{fig-1} (1)-(4) reveal the strength of structure property relationship between $SO_2$ and acentric factor, entropy, SHar and HNar, respectively.

By corresponding experimental values \cite{DTR21}, we get the correlations of $SO_2$ with some existing indices like Sombor index(SO), first ($M_1$) and second ($M_2$) Zagreb indices, forgotten topological index (F), connectivity index (R), sum connectivity index (SCI), symmetric division degree index (SDD) and neighborhood Zagreb index ($M_N$), which are listed in Table \ref{tb-3}.

\begin{table}
\begin{center}
	\caption{ The square of correlation coefficient of $SO_2$ with some existing indices} \label{tb-3}
	\begin{tabular}{c|cccccccc}
		\hline
&SO& $M_1$& $M_2$& F& R& SCI& SDD& $M_N$\\ \hline
$SO_2$&0.918&0.9201&0.8679&0.9022&0.8969&0.9150&0.8820&0.9123\\
\hline
	\end{tabular}
\end{center}
\small{Abbreviations: SO, Sombor index; SCI, sum connectivity index; SDD, symmetric division degree index.}
\end{table}

All these facts above show that the the second Sombor index is useful in predicting physicochemical properties.

\section{Extremal value of the second Sombor index of (Molecular) trees}
Let $G=(V,E)$ is a graph with order $n$ and size $m$. By the definition (\ref{eq-0}) of the second Sombor index of a graph, we have

\begin{equation}\label{eqt-1}
\begin{split}
SO_2(G)=\sum_{uv\in E}\frac{\frac{max\{d^2_G(u),d^2_G(v)\}}{min\{d^2_G(u),d^2_G(v)\}}-1}{\frac{max\{d^2_G(u),d^2_G(v)\}}{min\{d^2_G(u),d^2_G(v)\}}+1}=m-\sum_{uv\in E}\frac{2}{\frac{max\{d^2_G(u),d^2_G(v)\}}{min\{d^2_G(u),d^2_G(v)\}}+1}.
\end{split}
\end{equation}
Note that $1 \leq \frac{max\{d^2_G(u),d^2_G(v)\}}{min\{d^2_G(u),d^2_G(v)\}}\leq \frac{\Delta^2}{\delta^2}$, where $\delta$ and $\Delta$ are the minimum and maximum degree, respectively. From (\ref{eqt-1}), we straightforwardly obtain

\begin{thm} Let $G$ be a graph with $m$ edges and the maximal degree $\Delta$ and the minimal degree $\delta$, then
$$0\leq SO_2(G)\leq m\frac{\Delta^2-\delta^2}{\Delta^2+\delta^2}$$
the left equality holds if and only if every component of $G$ is regular and the right equality holds if and only if $\frac{max\{d_G(u),d_G(v)\}}{min\{d_G(u),d_G(v)\}}=\frac{\Delta}{\delta}$ for any edge $uv\in E$.
\end{thm}

\begin{thm}\label{thm-2} Let $T$ be a tree with $n>1$ vertices, then
$$\frac{6}{5}\leq SO_2(T)\leq \frac{(n^2-2n)(n-1)}{n^2-2n+2}$$
the left equality holds if and only if $T$ is the path $P_n$ and the right equality holds if and only if $T$ is the star $S_n$.
\end{thm}
\zm From (\ref{eqt-1}), we have
\begin{equation}
\begin{split}
SO_2(T)=n-1-\sum_{uv\in E}\frac{2}{\frac{max\{d^2_T(u),d^2_T(v)\}}{min\{d^2_T(u),d^2_T(v)\}}+1}.\nonumber
\end{split}
\end{equation}
And $\frac{max\{d_T(u),d_T(v)\}}{min\{d_T(u),d_T(v)\}}\leq n-1$, so
$$SO_2(T)\leq \frac{(n^2-2n)(n-1)}{n^2-2n+2}$$
with equality if and only if $T\cong S_n$.

On the other hand, we know that $f(x)=\frac{x^2-1}{x^2+1}$ is increasing for $x>0$. Let $E_1$ be the set of pendant edges in $T$ and $E_2=E-E_1$, then $|E_1|\geq 2$, $\frac{max\{d_T(u),d_T(v)\}}{min\{d_T(u),d_T(v)\}}\geq 2$ for $uv\in E_1$ and $\frac{max\{d_T(u),d_T(v)\}}{min\{d_T(u),d_T(v)\}}\geq 1$ for $uv\in E_2$. By (\ref{eqt-1}), we have
\begin{equation}
\begin{split}
SO_2(T)&=n-1-\sum_{uv\in E_1}\frac{2}{\frac{max\{d^2_T(u),d^2_T(v)\}}{min\{d^2_T(u),d^2_T(v)\}}+1}-\sum_{uv\in E_2}\frac{2}{\frac{max\{d^2_T(u),d^2_T(v)\}}{min\{d^2_T(u),d^2_T(v)\}}+1}\\
&\geq n-1-\frac{2}{5}|E_1|-|E_2|=\frac{3}{5}|E_1|\geq\frac{6}{5}.\nonumber
\end{split}
\end{equation}
with equality holds if and only if $T\cong P_n$.
\ezm

Let $\mathcal{C}T_n$ be the set of molecular trees with $n$ vertices. Let $m_{ij}$ be the number of edges with end vertices of degree $i$ and $j$ in $T$.

\begin{figure}[th]
	\centering
\includegraphics[width=0.8\textwidth]{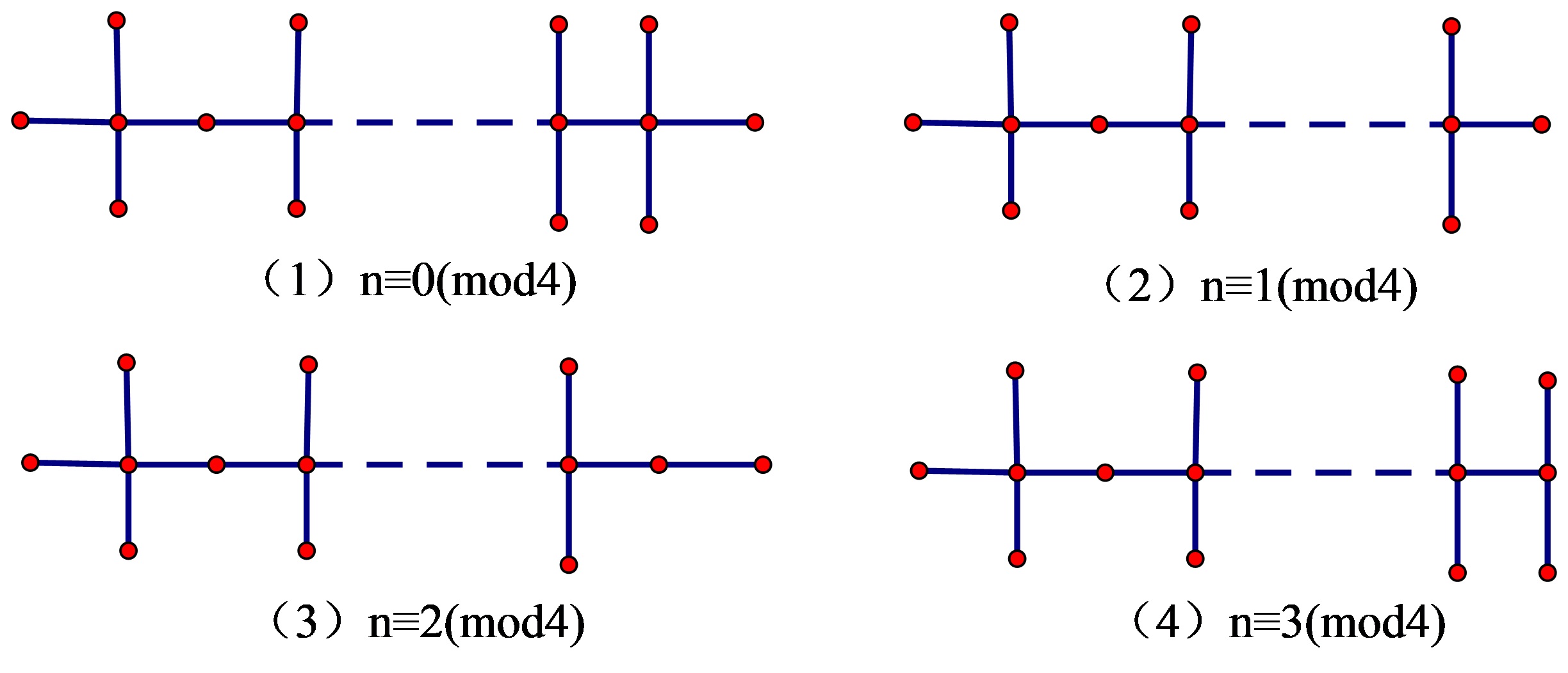}

	\caption{Four four types molecular trees with $n$ vertices }\label{fig-2}
\end{figure}

Denote by $\mathcal{T}_0$ the set of molecular trees $T$ with $n\equiv 0(mod 4)$ vertices, where $T$ has no vertex with degree 3, every vertex of degree 2 is adjacent to two vertices of degree 4 in $T$ and there is exactly a pair of vertices with degree 4 adjacent to each other, i.e., $T$ is a tree on $n$ vertices with $m_{14}=\frac{n+4}{2}$, $m_{24}=\frac{n-8}{2}$, $m_{44}=1$, $m_{12}= m_{13}= m_{22}= m_{23}= m_{33}= m_{34}=0$ (An example is shown in Figure \ref{fig-2} (1)).

Denote by $\mathcal{T}_1$ the set of molecular trees $T$ with $n\equiv 1(mod 4)$ vertices, where $T$ has no vertex with degree 3, every vertex of degree 2 is adjacent to two vertices of degree 4 in $T$ and no two vertices of degree 4 are mutually adjacent, i.e., $T$ is a tree with $m_{14}=\frac{n+3}{2}$, $m_{24}=\frac{n-5}{2}$, $m_{12}= m_{13}= m_{22}= m_{23}= m_{33}= m_{34}=m_{44}=0$ (An example is shown in Figure \ref{fig-2} (2)).

Denote by $\mathcal{T}_2$ the set of molecular trees $T$ with $n\equiv 2(mod 4)$ vertices, where $T$ has no vertex with degree 3, exactly one of vertices of degree 2 is adjacent to a vertex of degree 4 and a vertex of degree 1, and other vertices of degree 2 are adjacent to two vertices of degree 4, i.e., $T$ is a tree with $m_{14}=\frac{n}{2}$, $m_{24}=\frac{n-4}{2}$, $m_{12}=1$, $m_{13}= m_{22}= m_{23}= m_{33}= m_{34}=m_{44}=0$ (An example is shown in Figure \ref{fig-2} (3)).

Denote by $\mathcal{T}_3$ the set of molecular trees $T$ with $n\equiv 3(mod 4)$ vertices, where $T$ has exactly one vertex of degree 3, which is adjacent to two vertices of degree 1 and one vertex of degree 4 and every vertex of degree 2 is adjacent to two vertices of degree 4, i.e., $T$ is a tree with $m_{14}=\frac{n-1}{2}$, $m_{24}=\frac{n-7}{2}$, $m_{13}=2$, $m_{34}=1$, $m_{12}= m_{22}= m_{23}= m_{33}=m_{44}=0$ (An example is shown in Figure \ref{fig-2} (4)).

\begin{thm}\label{thm-3} Let $T\in \mathcal{C}T_n$ with $n\geq 5$, then
$$\frac{6}{5}\leq SO_2(T)\leq \left\{
\begin{array}{ll}
\frac{126n-108}{170}&n\equiv 0 (mod 4)\\
&\\
\frac{126n-30}{170}&n\equiv 1 (mod 4)\\
&\\
\frac{126n-102}{170}&n\equiv 2 (mod 4)\\
&\\
\frac{315n-281}{425}&n\equiv 3 (mod 4)
\end{array}
\right.
$$
the left equality holds if and only if $T\cong P_n$ and the right equality holds if and only if $T\in \mathcal{T}_i$ for $n\equiv i (mod 4)$ ($i=0,1,2,3$).
\end{thm}
\zm By Theorem \ref{thm-2}, we have $SO_2(T)\geq \frac{6}{5}$ with equality if and only if $T\cong P_n$.

From the definition (\ref{eq-0}) of the second Sombor index, we have
\begin{equation}\label{eqt-2}\begin{split} SO_2(T)=\frac{3}{5}m_{12}+\frac{4}{5}m_{13}+\frac{15}{17}m_{14}+\frac{5}{13}m_{23}+\frac{3}{5}m_{24}+\frac{7}{25}m_{34}.
\end{split}\end{equation}

Let $n_i$ be the number of vertices of degree $i$ in $T$, $i\in \{1,2,3,4\}$. We can get the following system of six linear equations which are satisfied by all molecular trees
\begin{equation}\left\{
\begin{split}
n_1+n_2+n_3+n_4&=n\\
n_1+2n_2+3n_3+4n_4&=2n-2\\
m_{12}+m_{13}+m_{14}&=n_1\\
m_{12}+2m_{22}+m_{23}+m_{24}&=2n_2\\
m_{13}+m_{23}+2m_{33}+m_{34}&=3n_3\\
m_{14}+m_{24}+m_{34}+2m_{44}&=4n_4
\end{split}\right.
\end{equation}
Solving this system with unknowns $m_{14}, m_{24}, n_1, n_2, n_3$ and $n_4$, we can obtain
\begin{equation}\label{eqt-5}
\left\{
\begin{split}
m_{14}=&\frac{n+3}{2}-\frac{3m_{12}}{2}-\frac{7m_{13}}{6}-\frac{m_{22}}{2}-\frac{m_{23}}{6}+\frac{m_{33}}{6}+\frac{m_{34}}{3}+\frac{m_{44}}{2}\\
m_{24}=&\frac{n-5}{2}+\frac{m_{12}}{2}+\frac{m_{13}}{6}-\frac{m_{22}}{2}-\frac{5m_{23}}{6}-\frac{7m_{33}}{6}-\frac{4m_{34}}{3}-\frac{3m_{44}}{2}\\
n_1=&\frac{n+3}{2}-\frac{m_{12}}{2}-\frac{m_{13}}{6}-\frac{m_{22}}{2}-\frac{m_{23}}{6}+\frac{m_{33}}{6}+\frac{m_{34}}{3}+\frac{m_{44}}{2}\\
n_2=&\frac{n-5}{4}+\frac{3m_{12}}{4}+\frac{7m_{13}}{12}+\frac{3m_{22}}{4}+\frac{m_{23}}{12}-\frac{7m_{33}}{12}-\frac{2m_{34}}{3}-\frac{3m_{44}}{4}\\
n_3=&\frac{m_{13}}{3}+\frac{m_{23}}{3}+\frac{2m_{33}}{3}+\frac{m_{34}}{3}\\
n_4=&\frac{n-1}{4}-\frac{m_{12}}{4}-\frac{m_{13}}{4}-\frac{m_{22}}{4}-\frac{m_{23}}{4}-\frac{m_{33}}{4}+\frac{m_{44}}{4}
\end{split}
\right.
\end{equation}
Replacing $m_{14}$ and $m_{24}$ in (\ref{eqt-2}) by (\ref{eqt-5}), we have
\begin{equation}\label{eqt-3}
\begin{split}
SO_2(T)=\frac{126n-30}{170}-\frac{36m_{12}}{85}-\frac{11 m_{13}}{85}-\frac{64 m_{22}}{85}-\frac{58 m_{23}}{221}-\frac{47 m_{33}}{85}-\frac{96 m_{34}}{425}-\frac{39 m_{44}}{85}
\end{split}
\end{equation}
which is maximal for a fixed number of vertices when the values $m_{12}, m_{13}, m_{22}, m_{23}, m_{33}$, $m_{34}$, and $m_{44}$ are equal to zero.
However, in the case of molecular trees with $n$ vertices, the condition
\begin{equation}\label{eqt-4}
\begin{split}
m_{12}= m_{13}= m_{22}= m_{23}= m_{33}= m_{34}=m_{44}=0
\end{split}
\end{equation}
can be satisfied only if $n\equiv 1(mod 4)$.

Any molecular tree satisfying (\ref{eqt-4}) has no vertices of degree 3, all its vertices of degree 2 are adjacent to two vertices of degree 4, and no two vertices of degree 4 are mutually adjacent (See (2) in Figure \ref{fig-2}).

Hence, if $n\equiv 1(mod 4)$, then for any molecular tree with $n$ vertices,
$$SO_2(T)\leq \frac{126n-35}{170}$$
with equality if and only if $T \in \mathcal{T}_1$.

If $n\not\equiv 1(mod 4)$, then (\ref{eqt-4}) cannot be satisfied by any molecular tree on $n$ vertices. In order to find the molecular trees with the maximal $SO_2$-value, we have to find the values of the parameters $m_{12}, m_{13}, m_{22}, m_{23}, m_{33}, m_{34}$, and $m_{44}$ as close to zero as possible compatible to the existence of a molecular tree, i.e., for which the right-hand sides of (\ref{eqt-5}) are integers, and for which a graph exists and we have that $m_{13}+m_{23}+2m_{33}+m_{34}$ has to be a multiple of 3 from $n_3=\frac{m_{13}}{3}+\frac{m_{23}}{3}+\frac{2m_{33}}{3}+\frac{m_{34}}{3}$.

By (\ref{eqt-3}), we know that there is must be $n_3=0$ or $n_3=1$ for $n\not\equiv 1(mod 4)$ if $T$ is a molecular tree with the maximal $SO_2(T)$-value.

{\bf Case 1.} If $n_3=0$, then $m_{13}=m_{23}=m_{33}=m_{34}=0$, and
\begin{equation}\label{eqt-6}
\begin{split}
SO_2(T)=\frac{126n-30}{170}-\frac{36m_{12}}{85}-\frac{64 m_{22}}{85}-\frac{39 m_{44}}{85}.
\end{split}
\end{equation}
To find the molecular tree(s) with the maximal $SO_2(T)$-value, we only need to consider $m_{12}+ m_{22}+ m_{44}=1$.

If $n\equiv 2(mod 4)$, there are two types of molecular trees such that $m_{12}+ m_{22}+ m_{44}=1$, i.e., $m_{12}=1, m_{22}= m_{44}=0$ or $m_{12}=0, m_{22}=1, m_{44}=0$. By simply computing and comparing, for any molecular tree $T$ with $n\equiv 1(mod 4)$ vertices,
$$SO_2(T)\leq \frac{126n-102}{170}$$
with equality if and only if $T \in \mathcal{T}_2$.

If $n\equiv 0(mod 4)$, there is only one type of molecular trees such that $m_{12}+ m_{22}+ m_{44}=1$, i.e., $m_{12}=m_{22}=0, m_{44}=1$. Then, for any molecular tree with $n\equiv 0(mod 4)$ vertices,
$$SO_2(T)\leq \frac{126n-108}{170}$$
with equality if and only if $T \in \mathcal{T}_0$.

If $n\equiv 3(mod 4)$, then there is no molecular trees such that $m_{12}+ m_{22}+ m_{44}=1$.

{\bf Case 2.} If $n_3=1$, then $m_{13}+m_{23}+2m_{33}+m_{34}=3$. To find the molecular trees with the maximal $SO_2(T)$-value, we only need to consider all possible choices of ($m_{13}$, $m_{23}$, $m_{33}$, $m_{34})$ such that $m_{13}+m_{23}+2m_{33}+m_{34}=3$ $m_{13}+m_{23}+2m_{33}+m_{34}=3$.

\begin{table}
\begin{center}
	\caption{ All ($m_{13}$, $m_{23}$, $m_{33}$, $m_{34})$ satisfy $m_{13}+m_{23}+2m_{33}+m_{34}=3$.} \label{tb-4}
	\begin{tabular}{c|c|c|c|c|c|c|c|c|c}
		\hline
$m_{13}$& $m_{23}$& $m_{33}$& $m_{34}$ & t &$m_{13}$& $m_{23}$& $m_{33}$& $m_{34}$ & t \\\hline
1&	1&	0&	1&				0.617715&
1&	0&	1&	0&				0.682341\\
0&	1&	1&	0&				0.815374&
0&	0&	1&	1&				0.778823\\
2&	1&	0&	0&				0.521233&
1&	2&	0&	0&				0.654266\\
1&	0&	0&	2&				0.581164&
2&	0&	0&	1&				0.484682\\
0&	2&	0&	1&				0.750748&
0&	1&	0&	2&				0.714197\\
\hline
	\end{tabular}
\end{center}
\small{where $t=\frac{11}{85}m_{13}+\frac{58}{221}m_{23}+\frac{47}{85}m_{33}+\frac{96}{425}m_{34}$.}
\end{table}

From Table \ref{tb-4} and (\ref{eqt-3}), a molecular tree on $n$ vertices and $n_3=1$ with the maximal $SO_2$-value must satisfy
\begin{equation}\label{eqt-9}
\begin{split}
m_{12}=m_{22}= m_{44}=m_{23}=m_{33}=0,m_{13}=2,m_{34}=1
\end{split}
\end{equation}
and it can be satisfied only if $n\equiv 3(mod 4)$. Then, for any molecular tree with $n\equiv 3(mod 4)$ vertices,
$$SO_2(T)\leq \frac{315n-281}{425}$$
with equality if and only if $T \in \mathcal{T}_3$.
\ezm

\bibliographystyle{elsarticle-num}
\bibliography{<your-bib-database>}

\begin{thebibliography}{00}

\bibitem{gut-1}I. Gutman, Geometric Approach to Degree-Based Topological Indices: Sombor Indices, MATCH Commun. Math. Comput. Chem. 86 (1) (2021) 11-16.

\bibitem{gut-2}I. Gutman, Sombor indices-back to geometry, Open Journal of Discrete Applied Mathematics, 5 (2) (2022) 1-5.

\bibitem{gut-3}I. Gutman, Degree-based topological indices. Croatica Chemica Acta, 86 (2013) 351-361.

\bibitem{Ku20}V. R. Kulli,. Graph indices. in: M. Pal, S. Samanta, A. Pal (Eds.), Handbook of Research of Advanced Applications
of Graph Theory in Modern Society. Hershey: Global, 2020, pp 66-91.

\bibitem{TC09}R. Todeschini, V. Consonni, Molecular Descriptors for Chemoinformatics.  Weinheim: Wiley¨CVCH, 2009.

\bibitem{DTR21}H. Deng, Z. Tang, R. Wu, Molecular trees with extremal values of Sombor indices, Int. J. Quantum Chem. 121 (2021) $\sharp$e26622.

\end{thebibliography}

\end{document}